\documentstyle[12pt,draft]{article}

 \newcounter{abceqn}

 \newcounter{abcfig}

\newcommand{\na}{\nabla}

\newcommand{\Om}{\Omega}

\newcommand{\pa}{\partial}

\newcommand{\Dl}{\Delta}

\newcommand{\ra}{\rightarrow}
\newcommand{\al}{\alpha}

\newcommand{\la}{\lambda}

\renewcommand{\theequation}{\thesection.\arabic{equation}}
\newcommand{\eqnsection}[1]{
	\section{#1}
	\setcounter{equation}{0}
	\renewcommand{\theequation}{\thesection.\arabic{equation}}
	\setcounter{figure}{0}
	\renewcommand{\thefigure}{\thesection.\arabic{figure}}
	\setcounter{remark}{0}
	\renewcommand{\theremark}{\thesection.\arabic{remark}}
	\setcounter{theorem}{0}
	\renewcommand{\thetheorem}{\thesection.\arabic{theorem}}
	\setcounter{lemma}{0}
	\renewcommand{\thelemma}{\thesection.\arabic{lemma}}
}

\title{{\bf Lax Pairs and Darboux Transformations for Euler Equations}}

\author{ \\ \\ 
Yanguang (Charles)\ \ Li  \thanks{This work is supported 
by the Guggenheim Fellowship. MSC: 76 34  35  37 }
\\ Department of Mathematics 
\\ University of Missouri \\ 
Columbia, MO 65211, USA \\ e-mail: cli@math.missouri.edu \\ 
phone: 573-884-0622 \\ fax: 573-882-1869
\\ \\ Artyom V. Yurov \\ Theoretical Physics Department \\
Kaliningrad State University, Russia \\ e-mail: yurov@freemail.ru}

\date{\today}

\renewcommand{\theequation}{\thesection.\arabic{equation}}

\begin{document}
\bibliographystyle{plain}
\maketitle
\newpage
\begin{abstract}    
In this article, we will report the recent developments on 
Lax pairs and Darboux transformations for Euler equations 
of inviscid fluids.
\end{abstract}

\newtheorem{lemma}{Lemma}
\newtheorem{theorem}{Theorem}
\newtheorem{corollary}{Corollary}
\newtheorem{remark}{Remark}
\newtheorem{definition}{Definition}
\newtheorem{proposition}{Proposition}
\newtheorem{assumption}{Assumption}

\tableofcontents



\eqnsection{Introduction}

The governing equations for the incompressible viscous fluid flow are the 
Navier-Stokes equations. Turbulence occurs in the regime of high Reynolds 
number. By formally setting the Reynolds number equal to infinity, the 
Navier-Stokes equations reduce to the Euler equations of incompressible 
inviscid fluid flow. One may view the Navier-Stokes equations with large 
Reynolds number as a singular perturbation of the Euler equations. 

Results of T. Kato show that 2D Navier-Stokes equations are globally 
well-posed in $C^0([0, \infty); H^s(R^2)), \ s>2$, and for any 
$0 < T < \infty$, the mild solutions of the 2D Navier-Stokes equations
approach those of the 2D Euler equations in $C^0([0, T]; H^s(R^2))$ 
\cite{Kat86}. 3D Navier-Stokes equations are locally well-posed in 
$C^0([0, \tau]; H^s(R^3)), \ s>5/2$, and the mild solutions of the 3D 
Navier-Stokes equations approach those of the 3D Euler equations in 
$C^0([0, \tau]; H^s(R^3))$, where $\tau$ depends on the norms of the initial 
data and the external force \cite{Kat72} \cite{Kat75}. Extensive studies on 
the inviscid limit have been carried by J. Wu et al. \cite{Wu96} 
\cite{CW96} \cite{Wu98} \cite{BW99}. There is no doubt that mathematical 
study on Navier-Stokes (Euler) equations is one of the most important 
mathematical problems. In fact, Clay Mathematics Institute has posted the 
global well-posedness of 3D Navier-Stokes equations as one of the one 
million dollars problems.

V. Arnold \cite{Arn66} realized that 2D Euler equations are a Hamiltonian 
system. Extensive studies on the symplectic structures of 2D Euler equations
have been carried by J. Marsden, T. Ratiu et al. \cite{Mar92}. S. Friedlander 
and M. Vishik \cite{FV90} \cite{VF93} found a Lax pair for Euler equations 
written in the Lagrangian coordinates.

Recently, the author \cite{Li00c} \cite{Li00d} found a Lax pair for 2D Euler 
equations written in the Eulerian coordinates. S. Childress \cite{Chi00}
found a Lax pair for 3D Euler equations. In this article, 
a Darboux transformation for 2D Euler equations and their Lax pair is 
found. Darboux transformations are powerful tools in generating explicit 
representations for figure eight structures in the phase spaces of 
Hamiltonian partial differential equations \cite{Li00a}. Of course,
Darboux transformations were traditionally used for generating 
multi-soliton solutions to soliton equations \cite{MS91}. 

Understanding the structures of solutions to Euler equations is of 
fundamental interest. Of particular interest is the question whether 
or not 3D Euler equations have finite time blow up solutions. T. Beale,
T. Kato, and A. Majda derived a necessary condition \cite{BKM84}.

Our hope is that the Lax pair and the Darboux transformation can be 
useful in investigating finite time blow up solutions of 3D Euler equations
and in establishing
global well-posedness of 3D Navier-Stokes equations. 

The rest of this article is organized as follows: In section 2, we will 
present a formulation of Euler equations. In section 3, we
will present a recent result on a Lax pair for 2D Euler equations. 
In section 4, we will present a recent result on a Lax pair for 3D 
Euler equations by Steve Childress. In section 5, we will present a 
recent result on 
a Darboux transformation for 2D Euler equations and their Lax pair.
In section 6, a short conclusion is presented.

\eqnsection{Formulations of Euler Equtions}

The three-dimensional incompressible Euler equation can be 
written in vorticity form,
\begin{equation}
\pa_t \Om + (u \cdot \na) \Om - (\Om \cdot \na) u = 0 \ ,
\label{3deuler}
\end{equation}
where $u = (u_1, u_2, u_3)$ is the velocity, $\Om = (\Om_1, \Om_2, \Om_3)$
is the vorticity, $\na = (\pa_x, \pa_y, \pa_z)$, 
$\Om = \na \times u$, and $\na \cdot u = 0$. $u$ can be 
represented by $\Om$ for example through Biot-Savart law.

When the fluid flow is of two dimensional, i.e. 
\[
u=(u_1(t,x,y), u_2(t,x,y), 0)\ ,
\]
the vorticity $\Om$ is of the form $\Om = (0, 0, \Om_3(t,x,y))$. 
Dropping the subscript $3$ of $\Om_3$, the 2D Euler equation can be 
written in the form,
\begin{equation}
\pa_t \Om + \{ \Psi, \Om \} = 0 \ ,
\label{euler}
\end{equation}
where the bracket $\{\ ,\ \}$ is defined as
\[
\{ f, g\} = (\pa_x f) (\pa_y g) - (\pa_y f) (\pa_x g) \ ,
\]
where $\Psi$ is the stream function given by,
\[
u=- \pa_y \Psi \ ,\ \ \ v=\pa_x \Psi \ ,
\]
and the relation between vorticity $\Om$ and stream 
function $\Psi$ is,
\[
\Om =\pa_x v - \pa_y u =\Dl \Psi \ .
\]

\eqnsection{A Lax Pair for 2D Euler Equation}

\begin{theorem}[Li, \cite{Li00c}]
The Lax pair of the 2D Euler equation (\ref{euler}) is given as
\begin{equation}
\left \{ \begin{array}{l} 
L \varphi = \la \varphi \ ,
\\
\pa_t \varphi + A \varphi = 0 \ ,
\end{array} \right.
\label{laxpair}
\end{equation}
where
\[
L \varphi = \{ \Om, \varphi \}\ , \ \ \ A \varphi = \{ \Psi, \varphi \}\ ,
\]
and $\la$ is a complex constant, and $\varphi$ is a complex-valued function.
\label{2dlp}
\end{theorem}
This theorem was announced in \cite{Li00c} without proof. Here we furnish 
a detailed proof.

Proof: In order to solve the over-determined system (\ref{laxpair}), one 
needs a compatibility condition. To derive a compatibility condition, we 
first take a $\pa_t$ to the first equation,
\[
(\pa_t L) \varphi + L (\pa_t \varphi) = \la (\pa_t \varphi) 
= - \la A \varphi = - A (\la \varphi) = -AL \varphi \ ,
\]
i.e.
\[
(\pa_t L) \varphi - LA \varphi = -AL \varphi \ ,
\]
and finally
\[
\bigg [ (\pa_t L) +[A, L] \bigg ] \varphi = 0\ ,
\]
where $[A, L] = AL - LA$ is the commutator. 
\[
(\pa_t L)\varphi =  \{ (\pa_t \Om), \varphi \} \ , 
\]
\begin{eqnarray}
[A, L]\varphi &=& \{ \Psi, \{ \Om, \varphi \} \} -
\{ \Om, \{ \Psi, \varphi \} \}  \nonumber \\
&=& \bigg [ \{ \Psi, \{ \Om, \varphi \} \} +
\{ \Om, \{ \varphi, \Psi \} \} \nonumber \\
& & + \{ \varphi, \{ \Psi, \Om \} \} \bigg ] - 
\{ \varphi, \{ \Psi, \Om \} \}\ , \nonumber
\end{eqnarray} 
where by Jacobi identity, $[ \ \cdot \  ] = 0$. Thus
\[
[A, L] \varphi = \{ \{ \Psi, \Om \}, \varphi \}\ .
\]
Finally
\[
\bigg [ (\pa_t L) +[A, L] \bigg ] \varphi = 
\{ \pa_t \Om + \{ \Psi, \Om \}, \varphi \} = 0\ .
\]
Therefore, the 2D Euler equation
\[
\pa_t \Om + \{ \Psi, \Om \} = 0
\]
is a compatibility condition. $\Box$

\begin{remark}
A Darboux transformation for the above Lax pair will be presented 
in section 5. Recently, Darboux transformations have been used in 
constructing explicit representations of homoclinic structures \cite{Li00a}.
\end{remark}

\eqnsection{A Lax Pair for 3D Euler Equation}

\begin{theorem}[Childress, \cite{Chi00}]
The Lax pair of the 3D Euler equation (\ref{3deuler}) is given as
\begin{equation}
\left \{ \begin{array}{l} 
L \varphi = \la \varphi \ ,
\\
\pa_t \varphi + A \varphi = 0 \ ,
\end{array} \right.
\label{3dlaxpair}
\end{equation}
where
\[
L \varphi = \Om \cdot \na \varphi - \varphi \cdot \na \Om \ , 
\ \ \ A \varphi = u \cdot \na \varphi - \varphi \cdot \na u \ , 
\]
$\la$ is a complex constant, and $\varphi = (\varphi_1, \varphi_2, 
\varphi_3)$ is a complex 3-vector valued function.
\end{theorem}
This theorem was obtained by Steve Childress \cite{Chi00}. For a proof, 
see the proof of Proposition \ref{prop}.

\begin{remark}
This Lax pair has the great potential that it can be very useful 
in investigating finite time blow up solutions of 3D Euler equation,
and in establishing global well-posedness of 3D Navier-Stokes equations.
\end{remark}
\begin{remark}
When one does 2D reduction to the Lax pair (\ref{3dlaxpair}),
one gets $L = 0$. Therefore, the Lax pair (\ref{3dlaxpair}) does not 
imply any Lax pair for 2D Euler equation.
\end{remark}
\begin{remark}
The most promising research following from this Lax pair 
should be along the direction of Darboux transformations, group 
symmetries etc.. In building the inverse scattering transform, it is 
crucial to have constant coefficient differential operators in 
(especially the spatial part of) the Lax pair. 
\end{remark}
\begin{proposition}[Li,\cite{Li00d}]
We consider the following Lax pair,
\begin{equation}
\left \{ \begin{array}{l} 
L \varphi = \la \varphi \ ,
\\
\pa_t \varphi + A \varphi = 0 \ ,
\end{array} \right.
\label{3vlaxpair}
\end{equation}
where
\[
L \varphi = \Om \cdot \na \varphi - \varphi \cdot \na \Om +D_1 \varphi\ , 
\ \ \ A \varphi = q \cdot \na \varphi - \varphi \cdot \na q +D_2 \varphi\ , 
\]
$\la$ is a complex constant, $\varphi = (\varphi_1, \varphi_2, 
\varphi_3)$ is a complex 3-vector valued function, $q=(q_1, q_2, q_3)$
is a real 3-vector valued function, $D_j = \al^{(j)} \cdot \na$, $(j=1,2)$,
and here $\al^{(j)}= (\al^{(j)}_1, \al^{(j)}_2, \al^{(j)}_3)$ are real constant
3-vectors. The following equation instead of the 3D Euler equation,
\begin{equation}
\pa_t \Om + (q \cdot \na) \Om - (\Om \cdot \na) q +D_2 \Om - D_1 q = 0 \ ,
\label{3veuler}
\end{equation}
is a compatibility condition of this Lax pair.
A specialization of (\ref{3veuler}) is the following system of equations,
\[
\left \{ \begin{array}{l} 
\pa_t \Om + (q \cdot \na) \Om - (\Om \cdot \na) q = 0 \ ,
\\ 
D_1 q = D_2 \Om\ .
\end{array} \right.
\]
\label{prop}
\end{proposition}
A chance in building an inverse scattering transform \cite{BC89} is that 
one can first build the inverse scattering transform 
for (\ref{3vlaxpair}), and then take the limits 
\[
\al^{(j)} \ra 0\ , \ (j=1,2)\ ,  \ \ \ q \ra u \ ,
\]
to get results for 3D Euler equation. 

Proof of Proposition \ref{prop}: Following the same argument as in the 
proof of Theorem \ref{2dlp}, one has the compatibility condition,
\[
\bigg [ (\pa_t L) + [A, L] \bigg ] \varphi = 0 \ .
\]
\[
(\pa_t L)\varphi = (\pa_t \Om) \cdot \na \varphi - \varphi
\cdot \na (\pa_t \Om)\ .
\]
\begin{eqnarray*}
[A, L]\varphi &=& q \cdot \na (\Om \cdot \na \varphi - 
\varphi \cdot \na \Om + D_1 \varphi) \\
& & - (\Om \cdot \na \varphi - \varphi \cdot \na \Om + D_1 \varphi)
\cdot \na q \\
& & + D_2 (\Om \cdot \na \varphi - \varphi \cdot \na \Om + D_1 \varphi) \\
& & - \Om \cdot \na ( q \cdot \na \varphi - \varphi \cdot \na q + 
D_2 \varphi ) \\
& & + ( q \cdot \na \varphi - \varphi \cdot \na q + D_2 \varphi )
\cdot \na \Om \\
& & - D_1 ( q \cdot \na \varphi - \varphi \cdot \na q + D_2 \varphi ) \\
&=& [(q \cdot \na)\Om]\cdot \na \varphi  
+ \underline{\Om \cdot [(q \cdot \na) \na \varphi ] \ \ }_{(1)}
- \underline{[(q \cdot \na)\varphi]\cdot \na \Om \ \ }_{(2)} \\
& & - \varphi \cdot [(q \cdot \na) \na \Om ]
+ \underline{(q \cdot \na)D_1 \varphi \ \ }_{(3)} 
- \underline{[(\Om \cdot \na)\varphi ]\cdot \na q \ \ }_{(4)}  \\
& & + [(\varphi \cdot \na)\Om ]\cdot \na q
- \underline{(D_1 \varphi)\cdot \na q \ \ }_{(5)} 
+ (D_2 \Om)\cdot \na \varphi \\
& & + \underline{(\Om \cdot \na)D_2 \varphi \ \ }_{(6)}
- \underline{(D_2 \varphi)\cdot \na \Om \ \ }_{(7)} 
- (\varphi \cdot \na)D_2 \Om \\
& & + \underline{D_2 D_1 \varphi \ \ }_{(8)}
- [(\Om \cdot \na)q]\cdot \na \varphi
- \underline{q \cdot [(\Om \cdot \na) \na \varphi ] \ \ }_{(1)} \\
& & + \underline{[(\Om \cdot \na)\varphi ]\cdot \na q \ \ }_{(4)}
+ \varphi \cdot [(\Om \cdot \na) \na q ]
- \underline{(\Om \cdot \na)D_2 \varphi \ \ }_{(6)} \\
& & + \underline{[(q \cdot \na)\varphi]\cdot \na \Om \ \ }_{(2)}
- [(\varphi \cdot \na) q]\cdot \na \Om
+ \underline{(D_2 \varphi)\cdot \na \Om \ \ }_{(7)} \\
& & - (D_1 q)\cdot \na \varphi
- \underline{(q \cdot \na)D_1 \varphi \ \ }_{(3)} 
+ \underline{(D_1 \varphi)\cdot \na q \ \ }_{(5)} \\
& & + (\varphi  \cdot \na)D_1 q
-\underline{D_1 D_2 \varphi \ \ }_{(8)}\ .
\end{eqnarray*}
All the terms $\underline{\ \ \ \ \ }_{(j)}, 1 \leq j \leq 8$, 
cancel each other. Thus, one has
\begin{eqnarray*}
[A, L]\varphi &=& \bigg \{ (q \cdot \na) \Om - (\Om \cdot \na)q
- D_1 q + D_2 \Om \bigg \} \cdot \na \varphi \\
& & - \varphi \cdot \na \bigg \{ (q \cdot \na) \Om - (\Om \cdot \na)q
- D_1 q + D_2 \Om \bigg \} \ .
\end{eqnarray*}
Therefore the following equation
\[
\pa_t \Om + (q \cdot \na) \Om - (\Om \cdot \na) q +D_2 \Om - D_1 q = 0 \ ,
\]
is a compatibility condition. Since this equation involves two vaiables 
$\Om$ and $q$, one can choose the specialization,
\[
\left \{ \begin{array}{l} 
\pa_t \Om + (q \cdot \na) \Om - (\Om \cdot \na) q = 0 \ ,
\\ 
D_1 q = D_2 \Om\ ,
\end{array} \right.
\]
to be a compatibility condition. $\Box$

\eqnsection{A Darboux Transformation for 2D Euler Equation}

Consider the Lax pair (\ref{laxpair}) at $\la =0$, i.e.
\begin{eqnarray}
& & \{ \Om, p \} = 0 \ , \label{d1} \\
& & \pa_t p + \{ \Psi, p \} = 0 \ , \label{d2} 
\end{eqnarray}
where we replaced the notation $\varphi$ by $p$.
\begin{theorem}
Let $f = f(t,x,y)$ be any fixed solution to the system 
(\ref{d1}, \ref{d2}), we define the Gauge transform $G_f$:
\begin{equation}
\tilde{p} = G_f p = \frac {1}{\Om_x}[p_x -(\pa_x \ln f)p]\ ,
\label{gauge}
\end{equation}
and the transforms of the potentials $\Om$ and $\Psi$:
\begin{equation}
\tilde{\Psi} = \Psi + F\ , \ \ \ \tilde{\Om} = \Om + \Dl F \ ,
\label{ptl}
\end{equation}
where $F$ is subject to the constraints
\begin{equation}
\{ \Om, \Dl F \} = 0 \ , \ \ \ \{ \Om +\Dl F, F \} = 0\ .
\label{constraint}
\end{equation}
Then $\tilde{p}$ solves the system (\ref{d1}, \ref{d2}) at 
$(\tilde{\Om}, \tilde{\Psi})$. Thus (\ref{gauge}) and 
(\ref{ptl}) form the Darboux transformation for the 2D 
Euler equation (\ref{euler}) and its Lax pair (\ref{d1}, \ref{d2}).
\label{dt}
\end{theorem}
\begin{remark}
For KdV equation and many other soliton equations, the 
Gauge transform is of the form \cite{MS91},
\[
\tilde{p} =  p_x -(\pa_x \ln f)p \ .
\]
In general, Gauge transform does not involve potentials.
For 2D Euler equation, a potential factor $\frac {1}{\Om_x}$
is needed. From (\ref{d1}), one has
\[
\frac{p_x}{\Om_x} = \frac{p_y}{\Om_y} \ .
\]
The Gauge transform (\ref{gauge}) can be rewritten as
\[
\tilde{p} = \frac{p_x}{\Om_x} - \frac{f_x}{\Om_x} \frac{p}{f}
=\frac{p_y}{\Om_y} - \frac{f_y}{\Om_y} \frac{p}{f}\ .
\]
The Lax pair (\ref{d1}, \ref{d2}) has a symmetry, i.e. it is 
invariant under the transform $(t,x,y) \ra (-t,y,x)$. The form 
of the Gauge transform (\ref{gauge}) resulted from the inclusion 
of the potential factor $\frac {1}{\Om_x}$, is consistent with 
this symmetry.
\end{remark}
\begin{proposition}
If $F$ satisfies the constraints
\begin{equation}
\{ \Om, \Dl F \} = 0 \ , \ \ \ \{ \Om, F \} = 0\ ,
\label{constraint1}
\end{equation}
or the constraints
\begin{equation}
\{ \Om, \Dl F \} = 0 \ , \ \ \ \{ \Dl F , F \} = 0\ ,
\label{constraint2}
\end{equation}
then $F$ satisfies the constraints (\ref{constraint}).
\end{proposition}

Proof: Notice that $\{ \Om, \Dl F \} = 0$ implies that
$\Dl F_x = \frac{\Om_x}{\Om_y}\Dl F_y$. Thus,
\begin{eqnarray*}
\{ \Om + \Dl F, F \} &=& (\Om_x + \frac{\Om_x}{\Om_y}\Dl F_y)F_y
- (\Om_y + \Dl F_y) F_x \\
&=& \frac{\Om_x}{\Om_y} (\Om_y + \Dl F_y) F_y -(\Om_y + \Dl F_y) F_x \\
&=& \frac{(\Om_y + \Dl F_y)}{\Om_y}\{ \Om, F \}\ .
\end{eqnarray*}
Similarly, one has
\[
\{ \Om + \Dl F, F \} =\frac{(\Om_y + \Dl F_y)}{\Dl F_y}\{ \Dl F, F \}\ .
\]
Thus the claim in the proposition is true. $\Box$
\begin{remark}
For soliton equations, $F$ can be represented in terms of $f$ \cite{MS91}.
\end{remark}

Proof of Theorem \ref{dt}: One notices that using (\ref{d1}), 
(\ref{d2}) can be rewritten as
\begin{equation}
p_t = \frac {p_x}{\Om_x}\{ \Om, \Psi \} \ .
\label{a0}
\end{equation}
To prove the theorem, one needs to check the two equations,
\begin{eqnarray}
& & \{ \tilde{\Om}, \tilde{p} \} = 0 \ , \label{ch1} \\
& & \tilde{p}_t = \frac {\tilde{p}_x}{\tilde{\Om}_x}
\{ \tilde{\Om}, \tilde{\Psi} \} \ . \label{ch2}
\end{eqnarray}
First we check (\ref{ch1}),
\[
\{ \tilde{\Om}, \tilde{p} \} = \{ \tilde{\Om}, \frac {p_x f - p f_x}
{\Om_x f} \} =0 \ ,
\]
leads to 
\begin{equation}
\frac {\tilde{\Om}_x \Om_y}{\Om_x^2 f^2}A - 
\frac {\tilde{\Om}_y \Om_x}{\Om_x^2 f^2}B = 0 \ ,
\label{abe}
\end{equation} 
where
\begin{eqnarray*}
A &=& \frac { (p_{xy}f + p_x f_y - p_y f_x - p f_{xy}) \Om_x f -
(p_x f - p f_x)(\Om_{xy} f + \Om_x f_y)}{\Om_y}\ , \\
B &=& \frac { (p_{xx}f - p f_{xx}) \Om_x f -
(p_x f - p f_x)(\Om_{xx} f + \Om_x f_x)}{\Om_x}\ .
\end{eqnarray*}
We will show that $A=B$. Using (\ref{d1}), one has
\[
\bigg (\frac {p_y}{\Om_y}\bigg )_x = \bigg (\frac {p_x}{\Om_x}\bigg )_x \ ,
\]
which leads to
\[
\frac {p_{xy}}{\Om_y} - \frac {p_y \Om_{xy}}{\Om_y^2}=
\frac {p_{xx}}{\Om_x} - \frac {p_x \Om_{xx}}{\Om_x^2} \ .
\]
Using (\ref{d1}) again,
\[
\frac {p_{xy}}{\Om_y} - \frac {p_x \Om_{xy}}{\Om_x\Om_y}=
\frac {p_{xx}}{\Om_x} - \frac {p_x \Om_{xx}}{\Om_x^2} \ .
\]
Similarly,
\[
\frac {f_{xy}}{\Om_y} - \frac {f_x \Om_{xy}}{\Om_x\Om_y}=
\frac {f_{xx}}{\Om_x} - \frac {f_x \Om_{xx}}{\Om_x^2} \ .
\]
Thus,
\[
\frac {p_{xy}f-pf_{xy}}{\Om_y}- \frac {(p_xf-pf_x)\Om_{xy}}{\Om_x\Om_y}
=\frac {p_{xx}f-pf_{xx}}{\Om_x} - \frac {(p_xf-pf_x)\Om_{xx}}{\Om_x^2}\ .
\]
That is,
\begin{equation}
\frac {(p_{xy}f-pf_{xy})\Om_x f}{\Om_y}=
\frac {(p_xf-pf_x)\Om_{xy}f}{\Om_y}+
\frac {(p_{xx}f-pf_{xx})\Om_x f}{\Om_x}
- \frac {(p_xf-pf_x)\Om_{xx}f}{\Om_x}\ .
\label{a1}
\end{equation}
Using (\ref{d1}) again, one has
\begin{equation}
\frac {p_xf_y-p_yf_x}{\Om_y}=p_x\frac {f_x}{\Om_x}
-\frac {p_x}{\Om_x}f_x = 0\ .
\label{a2}
\end{equation}
Using (\ref{a1}), (\ref{a2}), and (\ref{d1}), we have
\[
A = B = \frac {1}{\Om_x}\bigg [f^2 \Om_x p_{xx} - 
f(\Om_xf)_xp_x + [f_x (\Om_xf)_x - f_{xx}(\Om_xf)]p \bigg ] \ .
\]
Thus equation (\ref{abe}) becomes
\[
\{ \tilde{\Om}, \Om \}\frac {A}{\Om_x^2 f^2} = 0 \ .
\]
If we let
\begin{equation}
\{ \Om, \tilde{\Om} \} = \{ \Om, \Dl F\} = 0 \ ,
\label{act}
\end{equation}
then (\ref{abe}) is satisfied.

Next we check (\ref{ch2}). 
\[
\tilde{p}_t = \bigg (\frac {p_x f - p f_x}{\Om_x f}\bigg )_t\ .
\]
Using (\ref{a0}) for both $p$ and $f$, one gets
\begin{eqnarray}
\tilde{p}_t &=& \frac {1}{\Om_x^2 f^2} \bigg \{ \Om_x f 
\bigg [ f \bigg (\frac {p_x}{\Om_x}\{ \Om, \Psi \} \bigg )_x 
+ p_x \bigg (\frac {f_x}{\Om_x}\{ \Om, \Psi \} \bigg ) \nonumber \\
& & - f_x \bigg (\frac {p_x}{\Om_x}\{ \Om, \Psi \} \bigg )
- p \bigg (\frac {f_x}{\Om_x}\{ \Om, \Psi \} \bigg )_x \bigg ]\nonumber \\
& & - (p_x f - p f_x)\bigg [ f \{ \Om, \Psi \}_x + \Om_x
\frac {f_x}{\Om_x}\{ \Om, \Psi \} \bigg ] \bigg \} \nonumber \\
&=& (fp_{xx}-pf_{xx})\frac { \{ \Om, \Psi \} }{\Om_x^2 f} \nonumber \\
& & + (fp_{x}-pf_{x})\bigg [ \frac {1}{\Om_x f}
\bigg (\frac { \{ \Om, \Psi \} }{\Om_x} \bigg )_x - 
\frac {1}{\Om_x^2 f^2}\bigg (f\{ \Om, \Psi \}\bigg )_x \bigg ] \ . \label{b1}
\end{eqnarray}
The right hand side of (\ref{ch2}) is,
\begin{eqnarray}
\frac {\tilde{p}_x}{\tilde{\Om}_x}\{ \tilde{\Om}, \tilde{\Psi} \} 
&=& (fp_{xx}-pf_{xx})\frac {1}{\Om_x f}
\frac { \{ \Om + \Dl F, \Psi + F\} }{\Om_x + \Dl F_x}\nonumber \\
& & + (fp_{x}-pf_{x})\frac {- (\Om_x f)_x}{\Om_x^2 f^2}
\frac { \{ \Om + \Dl F, \Psi + F\} }{\Om_x + \Dl F_x} \ . \label{b2}
\end{eqnarray}
From (\ref{b1}) and (\ref{b2}), the coefficients of ``$fp_{xx}-pf_{xx}$'' 
leads to 
\begin{equation}
\frac {1}{\Om_x f} \frac { \{ \Om, \Psi \} }{\Om_x}
=\frac {1}{\Om_x f}\frac { \{ \Om + \Dl F, \Psi + F\} }{\Om_x + \Dl F_x}\ ,
\label{b3}
\end{equation}
and the coefficients of ``$fp_{x}-pf_{x}$'' leads to 
\begin{equation}
\frac {1}{\Om_x f}
\bigg (\frac { \{ \Om, \Psi \} }{\Om_x} \bigg )_x - 
\frac {1}{\Om_x^2 f^2}\bigg (f\{ \Om, \Psi \}\bigg )_x 
=\frac {- (\Om_x f)_x}{\Om_x^2 f^2}
\frac { \{ \Om + \Dl F, \Psi + F\} }{\Om_x + \Dl F_x} \ .
\label{b4}
\end{equation}
We will show that (\ref{b3}) implies (\ref{b4}). The left hand side of
(\ref{b4}) is
\begin{eqnarray*}
& & \frac {1}{\Om_x f}
\bigg (\frac { \{ \Om, \Psi \} }{\Om_x} \bigg )_x - 
\frac {1}{\Om_x^2 f^2}\bigg (\frac { \Om_xf\{ \Om, \Psi \} }
{ \Om_x }\bigg )_x \\
&=&\frac {1}{\Om_x f}
\bigg (\frac { \{ \Om, \Psi \} }{\Om_x} \bigg )_x -
\frac {1}{\Om_x^2 f^2}\Om_x f \bigg (\frac { \{ \Om, \Psi \} }
{ \Om_x }\bigg )_x \\
& & -\frac {1}{\Om_x^2 f^2}(\Om_x f)_x \frac { \{ \Om, \Psi \} }
{ \Om_x } = -\frac {1}{\Om_x^2 f^2}(\Om_x f)_x \frac { \{ \Om, \Psi \} }
{ \Om_x } \ .
\end{eqnarray*}
Thus (\ref{b4}) becomes 
\[
-\frac {1}{\Om_x^2 f^2}(\Om_x f)_x \frac { \{ \Om, \Psi \} }
{ \Om_x } =\frac {- (\Om_x f)_x}{\Om_x^2 f^2}
\frac { \{ \Om + \Dl F, \Psi + F\} }{\Om_x + \Dl F_x} \ .
\]
Therefore, (\ref{b3}) implies (\ref{b4}). From (\ref{b3}), we have
\[
(\Om_x + \Dl F_x)\{ \Om, \Psi \}  = \Om_x \{ \Om + \Dl F, \Psi + F\} \ ,
\]
which leads to
\begin{equation}
\Dl F_x \{ \Om, \Psi \}  = \Om_x [ \{ \Om, F\} +\{ \Dl F, F\} + 
\{ \Dl F, \Psi \} ] \ .
\label{b5}
\end{equation}
Notice that
\[
\Dl F_x \{ \Om, \Psi \}  -\Om_x \{ \Dl F, \Psi \}
= -\Psi_x \{ \Om, \Dl F\} = 0\ ,
\]
by (\ref{act}). Then (\ref{b5}) becomes
\begin{equation}
\Om_x [ \{ \Om, F\} +\{ \Dl F, F\} ] = 0 \ .
\label{b6}
\end{equation}
Thus if 
\[
\{ \Om + \Dl F, F\} = 0 \ ,
\]
(\ref{b6}) is fulfilled. $\Box$

\eqnsection{Conclusion}

We have briefly reported the most recent results on Lax pairs and a 
Darboux transformation for Euler equations of incompressible inviscid
fluid flows. The most promising researches following from these Lax 
pairs should be along the directions of Darboux transformations, 
group symmetries, etc. These Lax pairs have the great potentials that 
they can be very useful in investigating finite time blow up solutions 
of 3D Euler equations, and in establishing global well-posedness of 3D 
Navier-Stokes equations.

\bibliography{darboux}

\end{document}